\font\teneuf=eufm10 \font\seveneuf=eufm7 \font\fiveeuf=eufm5 \font\tenmsy=msbm10 \font\sevenmsy=msbm7 \font\fivemsy=msbm5 \font\tenmsx=msam10 \font\sixmsx=msam6 
\font\fivemsx=msam5 \textfont7=\teneuf \scriptfont7=\seveneuf \scriptscriptfont7=\fiveeuf \textfont8=\tenmsy \scriptfont8=\sevenmsy \scriptscriptfont8=\fivemsy \textfont9=\tenmsx \scriptfont9=\sixmsx \scriptscriptfont9=\fivemsx   
\noindent {\it This is the extended version of a paper that has appeared in the 
\vskip 1pt \noindent Mathematical Intelligencer. The final publication is available
\vskip 1pt
\noindent at Springer via http://dx.doi.org/10.1007/s00283-016-9644-3}
\vskip 20pt
\centerline{\bf Expect at most one billionth of a new Ferma{\it t} Prime!}  \vskip 20pt \centerline{Kent D. Boklan\footnote{\dag}{Department of Computer Science, Queens College, City University of New York boklan@cs.qc.cuny.edu} and John H. Conway\footnote{\ddag}{Gorenstein Visiting Professor, Department of Mathematics, Queens College, City University of New York and Department of Mathematics, Princeton University}}
\vskip 20pt \centerline{{\bf Abstract}: We provide compelling evidence that all Fermat primes were already known to Fermat.}   \vskip 20pt 
\noindent {\bf Prologue}
\vskip 10pt
\noindent What are the known Fermat primes? Hardy and Wright [HW] say that only four such primes are known, but this is incorrect since taking $F_n \ = \ 2^{2^n} + 1$, as they did, $F_0, F_1, F_2, F_3,$ and $F_4$ are prime. However, it's not clear that this is the definition that Fermat preferred.  Taking ``Ferma{\it t} prime" to 
mean ``prime of the form $2^n+1$", there are six known Ferma{\it t} primes, namely those for $n=0,1,2,4,8,16$. We shall pronounce the last letter of Ferma{\it t}'s name, as he did, when we include $2$ among the Ferma{\it t} primes, as he did. In this paper, we indicate this by italicizing that last letter, as we  already did.
\vskip 5pt \noindent Whichever definition we use, all the known such primes were already known to Fermat. Euler showed in 1732  that $2^{32}+1 \ = \ 4294967297= \ 641 \times 6700417$ is composite. In the endnotes to [HW], Hardy and Wright give a list of known composite Fermat numbers $F_n$, which is extended to 
$5 \leq n \leq 32$ with many other known composite values of $F_n$ (some with known factors, others merely known not prime); this is still the state of things. They go on to suggest that ``the number of primes $F_n$ is finite".  \vskip 5pt 
\noindent They then say ``This is what is suggested by considerations of probability ... The probability that a number $n$ is prime is at most $A / \log n$ and therefore the expected number of new Fermat primes is at most [a formula equivalent to]" $$A \sum {{1}\over{\log F_n }} < {{A}\over{\log 2}} \sum 2^{-n}  < 3A \ .$$
\vskip 8pt \noindent In this paper we produce compelling evidence for our thesis  (why should only Church have a thesis?):
\vskip 10pt \noindent {\bf Thesis}: {\it The only Ferma{\it t} primes are $2$ (according to taste), $3$, $5$, $17$, $257$ and $65537$.} 
\vskip 10pt \noindent As Hardy and Wright also say, their argument (notwithstanding its general lack of precision) assumes that there are no special reasons why a Fermat number should likely be a prime, while there are some. The most compelling ones are the result of Euler that every prime divisor of $F_n$ is congruent to $1$ 
modulo $2^{n+1}$ (we say ``$F_n$ is $2^{n+1}$-{\it full}"), and Lucas' 1891 strengthening of this to ``$F_n$ is $2^{n+2}$-{\it full}". A second reason is that the Fermat numbers are coprime, since $$F_{n+1} = F_0F_1F_2 ... F_n + 2 \ . $$
\vskip 10pt 
\noindent The Fermat number $F_n$ is either prime or not prime: the question of how to approximate the probability of primality for a general n is delicate. $F_n$ is not a generic odd but, if it were, according to the Prime Number Theorem (PNT), as there are about $n/{\log n}$ primes up to $n$, we would have $2/{\log n}$ as a naive first estimate for the primality of $F_n$. But the $F_n$ do have very special properties \ - \ and specific forms for prime factors. Since $F_n$ has no prime divisors that are $o(2^n)$, a simple (but loose) conditional probability estimate based on a lack of {\it small} divisors looks like $${{2}\over{\log n}} \prod_{p \le B} \left( 1 \ - \ {{1}\over{p}} \right)^{-1} $$ and by a classical result of Mertens, this is $$e^{\gamma} \ \ {{\log B}\over{\log n}} \ . $$ A tighter view is by way of the conditional probability that $F_n$ is prime given that all of its prime factors are congruent to $1 \ ({\rm mod} \ 2^{n+2})$; this is a basis for our discourse. 
\vskip 10pt
\noindent Prime divisors of the special form (as given by Lucas) seem to divide Fermat numbers with a greater likelihood than a generic prime might divide a generic number of the same size. Dubner and Keller [DK] comment that it ``appears that the probability of each prime [of the form] $k2^m+1$ [$k$ odd] dividing a standard Fermat number is $1/k$". If this were true, it may be contrasted with a prime $q$ of size $k2^m$ dividing $X$ of size $F_m$ as $${{1}\over{q}}\ =\ {{1}\over{k2^m}} \approx {{1}\over{k \log X}} \ \ , $$ in turn suggesting that the special form of the potential prime divisors of a Fermat number make them more likely to divide a Fermat number than a generic prime divisor. Our main result $(11)$ uses the most powerful conjectures about the distribution of primes in order to understand the limits of the method and
shows that even with the special form of possible prime factors, the probability that $F_n$ is prime is only at most a constant multiple larger than the naive estimate $2/{\log F_n}$, which demonstrates that the restriction of the possible divisors is nearly balanced out by the increased likelihood of each of these contender factors dividing $F_n$.
\vskip 12pt \noindent 
For some time every number-theorist has believed that the number of Fermat primes is finite. However, the Wikipedia page for Fermat Numbers  currently provides, but does not endorse, an argument in the 
reverse direction demonstrating the infinitude of the number of Fermat primes. (We invite the reader to find the flaw in the implied conditional probability argument that makes the estimate too large.)  
\vskip 12pt \noindent In the following section with details of our argument (the precision of which is the significant difference between ours and
previous heuristics), we take account of the main properties of Fermat numbers, concluding that the probability that $F_n$ is prime is approximately no larger than ${{4}\over{2^n}}$, which summed over $n \ge 33$ yields 
${{4}\over{2^{32}}} = {{1}\over{2^{30}}}$ as an upper bound for the expected number of new Ferma{\it t} primes, indeed less than one billionth. 
\vskip 20pt
\noindent {\bf An Argument for our Thesis}
\vskip 10pt \noindent Our objective is to estimate the probability that a number of size $x$ is a prime given only that it satisfies the various conditions that Fermat numbers are known to satisfy. 
\vskip 5pt \noindent  Let $\infty(x)$ denote a function that tends to infinity with $x$ and $\epsilon$ a fixed, arbitarily small positive number (that may vary in value from one use to the next). We write $\pi(I)$ for the number of primes in the interval $I$ and $\pi_{a(q)} (I)$ for the number of primes congruent to  $a ({\rm mod} \ q)$ 
therein.  For our argument, $I$ will be the interval $[x-r,x+r]$. We further set $\pi_{cond} (I)$ to be the number of primes in $I$ that satisfy a condition, ``$cond$". We write $\#_{cond}(I)$ for the number of integers in the interval $I$ satisfying the condition $cond$. And by $f(x) \sim g(x)$
we mean $f(x) = (1+o(1))g(x)$. 
\vskip 8pt \noindent We are concerned with the number of primes in the interval $[x-r,x+r]$ that satisfy various conditions. It is well-known that the number of primes in such an interval is approximately $2r/(\log x)$ provided $r$ is sufficiently large compared to $x$.  If the primes are restricted to 
be congruent to $a$ modulo $q$, then we expect a proportion $1/\phi(q)$ of this number if $r$ is large compared to $q$.
\vskip 5pt  \noindent We recall that the Fermat number $F_k = 2^{2^k}+1$ is $2^{k+2}$-{\it full}, that is all its prime divisors are congruent to $1 \ ({\rm mod} \ 2^{k+2})$. We desire to estimate the probability that a number $x$ (that satisfies some {\it cond}) in a certain interval $I$ is prime if we know that the number $x$ is $K$-{\it full},  in other 
words what we model as:  $${{\pi_{cond} (I)}\over{\#_{K-full} (I)}} \eqno(1)$$ where the length of the interval $I$ is taken as small as possible such that the expression in $(1)$ is meaningful.
\vskip 10pt \noindent {\bf I. Dealing with the $K$-{\it full} numbers}
\vskip 10pt \noindent If $x$ is large and $K$ is at least as large as $\log \log x$, it may be shown by a standard application of the Linear Sieve (see, for example, [HR]) that a first estimate for $\#_{K-full} [0,x]$ is
 $c \pi_{1 (K)} [0,x]$ where $c$ is a constant smaller than $2$\footnote{\dag}{It is worth noting  that a non-trivial lower bound on $K$ is needed, for if $K$ is bounded,  $\#_{K-full} [0,x]$ is genuinely of size $x$.}. The evaluation of the constant $c$ depends upon Mertens-type estimates for primes in arithmetic progressions.
(See [LZ].) This first estimate, applied to the interval $[x-r,x+r]$ where $x$ is large and fixed, still holds as the value of $r$ decreases from near $x$ to $(\log x)^{\delta}$, for some fixed $\delta$, and, meanwhile, the value of $c$ decreases to its limit of $1$. Numerical calculations kindly performed for us by Alex Ryba amply
confirm this and indeed suggest that the ratio $$ {{\#_{K-full} [x-r,x+r]}\over{\pi_{1 (K)} [x-r,x+r]}}$$  is much nearer to $1$ than $2$ for small $r$. So at the cost of at most a factor of $2$, and very likely less than $1.1$, we can replace``$K$-full number" by ``prime congruent to $1 ({\rm mod} \ K)$ in the denominator of 
$(1)$ and the resulting quotient is an upper bound for the probability in $(1)$.
\vskip 10pt \noindent {\bf II. Equidistribution if $r$ is large compared to $q$} 
\vskip 10pt \noindent So we are now faced with the question regarding primes in short intervals, which has been studied deeply. With no a priori information, the probability that a number of size $x$ be a prime is $(1/\log x)$ and this implies  that if $y$ is sufficiently large, an interval of length $y$ should contain about 
$y/(\log x)$ primes. Selberg [SE] showed, assuming the Riemann Hypthesis, that for almost all $x > 0$,  $$\pi(x+y) - \pi(x) \sim {{y}\over{\log x}}$$ holds provided that $$ y > \infty(x) \log^2 x   \eqno(2)$$ \vskip 2pt \noindent (where by {\it almost all} we mean in the sense of Lebesgue measure). 
\vskip 6pt \noindent  To generalize Selberg's result to arithmetic progressions $({\rm mod} \ q)$, as we need to do, we first require that the primes in $I$ be uniformly distributed among the $\phi(q)$ residue classes provided only that there are sufficiently many of them. That is, provided $r$ is large 
enough compared to $q$. This is solved by
\vskip 10pt \noindent {\bf The Equidistribution Lemma}: {\it If $B$ balls are independently distributed into $C$ cups then provided that  $B > \infty(C) C \log C$ one expects that for any fixed $\epsilon > 0$, the number of balls contained in any one cup is between $(1-\epsilon)(B/C)$ and $(1 + \epsilon)(B/C)$.}
\vskip 10pt \noindent (We thank Noam D. Elkies [NDE] for an argument establishing this Lemma - a result we have not been able to locate in the literature although, as Elkies remarks and we agree, it must be ``well-known".)
\vfil \eject
\vskip 6pt \noindent This Lemma shows that the primes in $I=[x-r,x+r]$ will be uniformly distributed (in the above sense) provided that $${{r}\over{\log x}} > \infty(x) \phi(q) \log \phi(q) \ .  $$ 
We now set $q=(\log x)^{\delta}$, with $\delta \ge 1$, so that the Lemma's requirement - and our requirement on $I$ - becomes $$r > (\log x)^{1 + \delta + \epsilon} \ . \eqno(3) $$ (Note that this implies that $r > \infty(x) \log^2 x$, so
Selberg's condition $(2)$ is automatically satisfied and our application of the Lemma with $r/(\log x)$  ``balls" is justified.) 
\vskip 10pt \noindent {\bf III. Uniformity if $r$ is large compared to $x$}
\vskip 10pt \noindent In order to establish a lower bound condition on $r$ (half the length of $I$) for the expected result on primes in arithmetic progression 
$$\pi_{1(q)}(I) \sim  {{2r}\over{\phi(q) \log x}} \eqno(4)$$ to hold we first set $$I(x,r,q,a) = \int_x^{2x} \left( \psi(y+h,q,a) - \psi(y,q,a) - {{h}\over{\phi(q)}} \right)^2
dy$$ where $$\psi(y,q,a) = \sum_{{{n \le x}\atop{n \equiv a (q)}}} \Lambda(n) \ , $$ and $\Lambda(n)$ is the von Mangoldt function (defined by $\Lambda(n)=\log p$ if $n$ is a power of the prime $p$ and $0$ otherwise). Then we define $$I(x,h,q) = \sum_{a({\rm mod} \ q)} I(x,h,q,a)$$ where 
``$a \ ({\rm mod} \ q)$" indicates that the sum is taken over all residue 
classes  $a \ ({\rm mod} \ q)$ with $(a,q)=1$. On the Generalized Riemann Hypothesis (GRH), Prachar [PR] showed that\footnote{\ddag}{In the other direction, assuming the GRH, for $x \ge 2$, $1 \le q \le h \le x$ and $q \le h \le (xq)^{1/3-\epsilon}$, Goldston and Yildirim [GY] have shown that $$I(x,h,q,a) \ge {{1}\over{2}} 
{{xh}\over{\phi(q)}} \log \left({{xq}\over{h^3}} \right) ( 1 - o(1)) \ . $$ Combining this with the Equidistribution Lemma, Prachar's result $(5)$ is at most $1-\epsilon$ logarithms away from being best possible.}  $$I(x,h,q) =  O \left(hx \log^2 (qx) \right)  \ . \eqno(5)$$ 
\vskip 6pt \noindent We wish in some way to identify the 
$a \ ({\rm mod} \ q)$ with $(a,q)=1$ for which for which the expected asymptotic relation $$\psi(y+h,q,a) - \ \psi(y,q,a) \sim  {{h}\over{\phi(q)}} \eqno(6)$$ fails to hold (for a positive proportion of $y$ in $[x,2x]$). Using $\sum ^*$ to indicate summation over such $a$, note that 
$$I(x,h,q) \ge \sum \ ^* \ \int_x^{2x} \left( \psi(y+h,q,a)  - \psi(y,q,a) - {{h}\over{\phi(q)}} \right)^2 dy$$ \noindent and so, from $(5)$, $$\sum \ ^* \ 1 = O  \left( {{\phi^2(q) \log^2 x}\over{h}} \right) \ . \eqno(7)$$ \noindent 
It follows immediately that if  $h > \infty(x) \phi(q) \log^2 x$, then the asymptotic estimate equivalent to $(6)$, namely $$\pi(y+h,q,a) - \pi(y,q,a) \sim {{h}\over{\phi(q) \log y}} \eqno(8)$$ \noindent holds for almost all $a$ mod $q$ for fixed $q$ for almost all $y$ in $[x,2x]$.  \vskip 5pt \noindent At this point we introduce what 
we call the \vskip 10pt \noindent {\bf Uniformity Conceit}: {\it If an assertion holds for almost all values of $a$ it is very probable that it holds for any particular $a$ unless there is a good reason why it should not.}  
\vskip 10pt \noindent We apply the Conceit with $h > \infty(x) \phi(q) \log^2 x$ to $(8)$ with $a=1$ to deduce that very probably $$\pi(x+h,q,1) - \pi(x,q,1) \sim {{h}\over{\phi(q) \log x}}$$ holds. Setting $h=2r$ and shifting $x$ to $x-r$ we find that very probably $(4)$ holds so long as the following condition 
(stricter than $(3)$which was needed for equidistribution across residue classes) is satisfied: $$r > (\log x)^{2+\delta + \epsilon} \ . \eqno(9)$$
\vskip 10pt \noindent {\bf IV. Towards the probability that $F_n$ is prime} 
\vskip 10pt \noindent So now we return to the probability that a number $N$, satisfying condition $cond$, in $I=[x-r,x+r]$ is a prime given that it is $K$-{\it full}. From the probability given in $(1)$, we take $x$ to be the Ferma{\it t} number $F_n$ with $K=2^{n+2}$. Then using our {\it 
replacement} of ``$K$-full" by ``prime congruent to $1 \ ({\rm mod} \ K)$", an upper bound for this probability is $${{\pi_{cond} (I)}\over{\pi_{1(2^{n+2})}(I)}} \eqno(10)$$ for $r$ as small as possible such that this is meaningful. 
\vskip 6pt \noindent The Fermat number $F_n$ is $2^{n+2}$-{\it full} so $F_n$ is in the arithmetic progression $1 \ ({\rm mod} \ 2^{n+2})$.  Of course, we know that $F_n \equiv 1 \ ({\rm mod} \ 2^{2^{n}})$ but we cannot use a value of $q=2^{2^n}$ because The Equidistribution Lemma requires that
$$\phi(q) \ = \  o \left( {{x}\over{\log^2 x}} \right) \ . $$ However, we do need and can use something larger than $2^{n+2}$ and can in fact use $q=2^{2n}$ so $F_n$ lies in the progression $1 \ ({\rm mod} \ 2^{2n})$; 
we take this as condition $cond$ in the numerator of $(10)$.\footnote{\dag}{Selecting a larger value of $q$ such as $2^{\alpha n}$ for $\alpha > 2$ leads to a smaller upper bound for the probability we seek, but increasing the value of $\alpha$ requires a larger $r$ which is discordant  with our model, making the interval 
around $x$ as small as possible.} \vskip 6pt
\noindent On all of our reasonable assumptions (including the GRH), we may evaluate this probability provided $r$ satisfies $(9)$\footnote{\ddag}{The observant reader will have noticed that we have not mentioned the second of Hardy and Wright's reasons that the Fermat numbers are more likely to be prime than other
numbers, namely that they are mutually coprime; but that observant reader will also have noticed that the Fermat numbers increase so rapidly that at no time in this Argument has this made any difference.}. 
For $q=(\log x)^{\delta} \sim 2^{2n}$ we have $\delta=2$ so, from $(9)$, the interval $I$ must be at least as  large as $2 (\log x)^{4+\epsilon}$ \ ; we select $r$ to be $(\log x)^{4+\epsilon}$. Thus, from $(4)$, in accord with our model, the probability that $F_n=x$ is prime is at most 
$${{\pi_{1 (2^{2n})} [x-r,x+r]}\over{\pi_{1 (2^{n+2})} [x-r,x+r]}} \sim {{{{1}\over{\phi(2^{2n})}} {{2r}\over{\log x}}}\over{{{1}\over{\phi(2^{n+2})}} {{2r}\over{\log x}}}} \ =  \ {{4}\over{2^n}} \  . \eqno(11) $$ For a new Fermat prime $F_n$, $n$ must be at least $33$, so the expected number of new Fermat primes is at
most $$ {{4}\over{2^{33}}} + {{4}\over{2^{34}}} + {{4}\over{2^{35}}} + \dots = {{4}\over{2^{32}}}$$ and since $2^{30}$ exceeds one billion, this is indeed less than one-billionth, justifying the title of our paper.
\vskip 20pt \noindent {\bf V.  Improvements?}
\vskip 10pt \noindent We - your authors - are firmly convinced that there will be no significant imrovement on our paper throughout all of future time.
This seems an extremely audacious statement, but the evidence we now present should make it extremely plausible. 
\vskip 6pt \noindent  Some small progress will undoubtedly be made.  We divide the numbers $F_n$ ($n \ge 33$) into two lists
\vskip 6pt (A) 36,37,38,39,42,43,... 
\vskip 6pt \noindent the numbers for which a prime factor of $F_n$ has been found and 
\vskip 6pt (B) 33,34,35,40,41,...
\vskip 6pt \noindent for which nothing is known. 
The list (A) has been formed by taking a prime $p$ of the shape $k 2^m +1$ where $k$ is sufficiently small and $m$ sufficiently large and
finding that the sequence $2,2^2,2^4,2^8...$ formed by repeated squaring ({\rm mod} $p$). If the number $2^{2^r}$ is conguent to $-1$
mod $p$ then one has shown that $F_r$ is divisble by $p$ so that $r$ should belong to list (A).
\vskip 10pt \noindent
This procedure will obviously be continued; so that numbers will gradually be promoted to list (A) from list (B). Your authors, who have not memorized
the two lists, will not regard this as significant progress.
\vskip 6pt \noindent There is a necessary and sufficient test for primality of $F_n$ but it involves so much computation that the largest $F_n$ that
has been proved composite this way is $F_{24}$ [CMP]. For this, the computation took the better part of a year of computer time in 1999. 
Now $F_{33}$ is $2^9$ times as long as $F_{24}$, it will may well fall to this method at least when quantum computation fulfills the
promises that have been made for it. If so, some of our readers may live to see a paper that updates our title to 
"Expect at most one trillionth ....". Will you regard this as a substantial improvement?
\vskip 20pt
\noindent {\bf Epilogue}
\vskip 10pt
\noindent The famous question of the infinitude of the Mersenne primes falls in the orbit of our approach; setting $M_p=2^p-1=x$ it is easy to find that $M_p$ is $p-full$. When it comes to Mersenne primes, naive (PNT) estimates suggest there are infinitely many because $\sum 1/p$ diverges; we do not quibble with this expectation. But the approach in the Argument is not directly amenable to the Mersenne prime question since we lack a {\it second condition}. Adding in such a condition into our model leads to our
\vskip 10pt
\noindent {\bf Conjecture:} There are finitely many Mersenne primes $M_p$ whose index $p$ is a twin prime.
\vskip 10pt
\noindent The evidence is equally strong for the
\vskip 10pt
\noindent {\bf Conjecture:} There are finitely many Mersenne primes $M_p$ whose index $p$ is a Sophie Germain prime (i.e. a prime $p$ where $2p+1$ is also prime).
\vskip 10pt
\noindent More generally, if we fix integers $a$ and $b$ where $(b,p)=1$ and consider primes of the form $M_p$ where both $p$ and $ap+b$ are primes, the Selberg Sieve (see [HR]) establishes an upper bound on the count of such primes up to $x$, $$\sum_{p < x} 1 \ = \ O \left({{x}\over{\log^2 x}} \right) \ .$$ We may apply the technique in the Argument, with $M_p$ taking on $a$ distinct residues modulo $ap+b$ so the Chinese Remainder Theorem determining $a$ values of $M_p$  modulo $p(ap+b)$, to produce an upper bound probability of size $1/(ap+b)$ for the probability that such an $M_p$ is prime. Summing over the special primes using partial summation yields a total expectation of $O(1)$ of these special Mersenne primes. So we are led to formulate the general
\vskip 10pt
\noindent {\bf Conjecture:} There are only finitely many Mersenne primes $2^p-1$ where $p$ is a prime and $ap+b$ is also prime (for some fixed integers $a$ and $b$ where $(b,p)=1$).   
\vskip 25pt
\noindent We thank an anonymous referee to an earlier version of this paper for several helpful suggestions. 
\vskip 20pt \centerline{\bf References} \vskip 15pt
\noindent [CMP] Crandall, R., Mayer, E. and Papadopoulos, J.: The Twenty-Fourth Fermat Number is Composite. {\it Math. of Computation} {\bf 72} (2002), 1555-1572 \vskip 10pt
\noindent [DK] Dubner, H. and Keller, W.: Factors of generalized Fermat numbers. {\it Math. of Computation} {\bf 64} (1995), 397-405
\vskip 10pt \noindent [GY] Goldston, D. and Yildirim, C.: On the Second Moment for primes in an arithmetic progression. {\it Acta Arith.} {\bf 100} (2001), 85-104
\vskip 10pt \noindent [HR] Halberstam, H. and Richert, H.-E.: {\it Sieve Methods}. New York, Academic Press, 1974 
\vskip 10pt \noindent [HW] Hardy, G.H. and Wright, E.M.: {\it An Introduction to the Theory of Numbers} [6th ed. (2008)] Oxford University Press
\vskip 10pt  \noindent [LZ] Languasco, A. and Zaccagnini, A.: On the constant in the Mertens product for arithmetic progressions II: Numerical Values. 
{\it Math. of Computation} {\bf 78} (2009), 315 - 326
\vskip 10pt \noindent [NDE] Private email communication to the first author, January 3, 2014
\vskip 10pt  \noindent [PR] Prachar, K.: Generalisation of a theorem of A. Selberg on primes in short intervals. In: {\it Topics in Number Theory, Colloquia Mathematica Societatis Ja\'nos Bolyai 13}. Debrecen, 1974, 267-280
\vskip 10pt  \noindent [SE] Selberg, A.: On the normal density of primes in small intervals, and the difference between consecutive primes. {\it Archiv for Mathematik og Naturvidenskap} {\bf 47} (1943), 87-105 
\end